\documentclass[oneside,notitlepage,12pt]{article}

\usepackage{amssymb}
\usepackage[leqno]{amsmath}
\usepackage{amsfonts}
\usepackage{amsopn}
\usepackage{amstext}
\usepackage{amsthm}
\usepackage{enumitem}
\usepackage[all]{xy}
\newdir{ >}{{}*!/-9pt/@{>}}

\usepackage[colorlinks, backref]{hyperref}

\usepackage{calrsfs}
\usepackage{fourier-orns}
\usepackage{hieroglf} 
\usepackage{clock} 

\newcommand{\rloe}{\preceq}

\providecommand{\cal}{\mathcal}
\renewcommand{\Bbb}{\mathbb}

\newenvironment{pf}{\begin{proof}}{\end{proof}}

\newcommand{\Dee}{{\cal{D}}}

\newcommand{\Tau}{{\cal{T}}}
\newcommand{\Yu}{{\cal{U}}}

\newcommand{\Nat}{{\Bbb{N}}}

\newcommand{\al}{\alpha}

\newcommand{\sig}{\sigma}

\renewcommand{\phi}{\varphi}
\renewcommand{\rho}{\varrho}

\newcommand{\unii}{{\mathbb I}}
\newcommand{\ntr}{{n\in\omega}}
\newcommand{\Ntr}{n\in{\Bbb{N}}}
\newcommand{\loe}{\leqslant}
\newcommand{\goe}{\geqslant}

\newcommand{\subs}{\subseteq}
\newcommand{\sups}{\supseteq}
\newcommand{\nnempty}{\ne\emptyset}

\newcommand{\ovr}{\overline}

\newcommand{\cl}{\operatorname{cl}}

\newcommand{\id}[1]{{\operatorname{i\!d}_{#1}}} 

\newcommand{\setof}[2]{\{#1\colon #2\}}

\newcommand{\sett}[2]{\{#1\}_{#2}}
\newcommand{\sn}[1]{\{#1\}} 
\newcommand{\pair}[2]{\langle #1, #2 \rangle} 
\newcommand{\map}[3]{#1\colon #2 \to #3} 
\newcommand{\img}[2]{#1[#2]} 

\newcommand{\fra}{Fra\"iss\'e}

\newcommand{\U}{\mathbb U}

\providecommand{\nat}{\omega}

\newcommand{\ciag}[1]{{\sett{{#1}_n}{\ntr}}}

\newcommand{\ciagi}[1]{\sig{#1}}

\newcommand{\fK}{{\mathfrak{K}}}
\newcommand{\fL}{{\mathfrak{L}}}

\newcommand{\fG}{{\mathfrak{G}}}

\newcommand{\cmp}{\circ} 

\newcommand{\separator}{\begin{center}***\end{center}}

\newcommand{\br}[1]{\partial#1}

\newcommand{\ciint}[1]{\left[{#1},\rightarrow\right)}

\newcommand{\icint}[1]{\left(\leftarrow,{#1}\right]}

\newcommand{\proto}[1]{{\mathbb S_\kappa}}

\newtheorem{tw}{Theorem}[section]
\newtheorem{wn}[tw]{Corollary}

\newtheorem{prop}[tw]{Proposition}

\theoremstyle{definition}
\newtheorem{df}[tw]{Definition}
\newtheorem{ex}[tw]{Example}

\theoremstyle{remark}

\newcommand{\BM}[1]{\operatorname{BM}\left(#1\right)}
\newcommand{\BMG}[2]{\BM {#1, #2} }

\title{Banach-Mazur game played in partially ordered sets}
\author{
{\sc Wies{\l}aw Kubi\'s}\footnote{Research supported by NCN grant DEC-2011/03/B/ST1/00419.}\\ \\
{\small Department of Mathematics}\\ {\small Jan Kochanowski University in Kielce, Poland}\\
{\small \textit{and}}\\
{\small Institute of Mathematics, Academy of Sciences of the Czech Republic}
}
\date{\today\ \clocktime}

\begin{document}

\maketitle

\begin{abstract}
We present a version of the Banach-Mazur game, where open sets are replaced by elements of a fixed partially ordered set.
We show how to apply it in the theory of \fra\ limits and beyond, obtaining simple proofs of universality of certain objects and classes.


\end{abstract}

\tableofcontents

\section{Introduction}

The Banach-Mazur game is usually played in a topological space, using its nonempty open subsets.
The idea is that two players alternately build a decreasing sequence of sets and the result of the game is its intersection.

In this note we develop a more abstract setting for this game.
Namely, the family of all nonempty open sets in a fixed topological space can be regarded as a partially ordered set.
In order to say who wins, one needs to distinguish a ``winning" family of countably generated ideals of this poset.
More precisely, one of the players wins if the ideal generated by the sequence resulted from a play belongs to our distinguished family.
It turns out that one can reformulate and extend known results in this new setting (Section~\ref{SecKlasik}).
As an application, we discuss the Banach-Mazur game played with finitely generated models taken from a \fra\ class, showing that one of the players has a winning strategy ``leading to" the \fra\ limit (Section~\ref{SecFLimswegb}).

\section{Preliminaries}

Here we put the relevant concepts, notions, definitions, and basic facts.

Given a topological space $X$, we denote by $\Tau^+(X)$ the collection of all nonempty open subsets of $X$.
Recall that a \emph{$\pi$-base} in $X$ is a family $\Yu \subs \Tau^+(X)$ such that for every $V \in \Tau^+(X)$ there is $U \in \Yu$ satisfying $U \subs V$.
In this note, topological spaces are not assumed to satisfy any separation axioms, except $T_0$ (that is, open sets should separate points).

We shall use some basic notions concerning partially ordered sets (briefly: posets).
Namely, given a poset $\pair P \loe$, a subset $D \subs P$ is \emph{cofinal} in $P$ if for every $p \in P$ there is $d \in D$ with $p \loe d$.
Note that a cofinal subset of $\pair {\Tau^+(X)}{\sups}$ is just a $\pi$-base of the topological space $X$.
An \emph{ideal} in $P$ is a set $I \subs P$ satisfying the following conditions:
\begin{enumerate}
\item[(I1)] $(\forall\; x,y \in I)(\exists\; z \in I)\;\; x \loe z \text{ and } y \loe z$;
\item[(I2)] $(\forall\; x \in I) \icint x \subs I$.
\end{enumerate}
An ideal $I$ is \emph{countably generated} if it has a countable cofinal subset.
Note that $I$ is countably generated if and only if there is an increasing sequence $\ciag a$ such that
$$I = \setof{x \in P}{(\exists\; \ntr)\;\; x \loe a_n}.$$
We shall denote by $\ciagi P$ the poset whose elements are all countably generated ideals of $P$ and the ordering is inclusion.
Note that $\ciagi P$ is a natural extension of $P$, namely, each $p\in P$ can be identified with $\icint p \in \ciagi P$.
Note also that every increasing sequence in $\ciagi P$ has the supremum in $\ciagi P$.
In fact, $\ciagi P$ can be called the ``sequential completion'' of $P$, because of the following reason. If $\map f P Q$ is order preserving and $Q$ is sequentially complete (that is, every increasing sequence in $Q$ has the supremum) then there is a unique extension $\map {\tilde f} {\ciagi P}{Q}$	
of $f$ to a sequentially continuous order preserving mapping, given by the formula $\tilde f(I) = \sup I$.
Recall that a mapping $g$ between sequentially complete posets is \emph{sequentially continuous} if $\sup_{\ntr}g(x_n) = g(\sup_{\ntr}x_n)$ for every increasing sequence $\ciag x$.

\subsection{Metric trees}

Recall that a \emph{tree} is a poset $\pair T \loe$ such that for every $t \in T$ the set
$$\setof{x \in T}{x < t}$$
is well-ordered. Its order type is the \emph{height} of $t$ in $T$.
The set of all elements of $T$ of a fixed height $\al$ is called the $\al$th \emph{level} of $T$.
The \emph{height} of a tree $T$ is the minimal ordinal $\delta$ such that level $\delta$ is empty, i.e. every element of $T$ has height $<\delta$.
Maximal elements of a tree are called \emph{leaves}.
We are interested mainly in trees of height $\omega$.
We shall call them \emph{metric trees}.
We use the adjective \emph{metric} in order to emphasize that such a tree naturally leads to a metrizable space.

Namely, given a tree $T$, let $\br T$ denote the set of all \emph{branches} of $T$, that is, all maximal chains in $T$.
A branch in $T$ may have a maximal element, called a \emph{leaf}.
Given $t \in T$, we set $t^+ = \setof{X \in \br T}{t \in X}$.
Then the family $\sett{t^+}{t \in T}$ is a basis of a topology on $\br T$, and the sets $t^+$ clopen (i.e., closed and open) with respect to this topology.
In particular, $\br T$ is zero-dimensional.
If $T$ is a metric tree then this topology is metrizable.
Indeed, given $X, Y \in \br T$, if $X \ne Y$ then we may define their distance $\rho(X,Y)$ to be $1/n$ where $n$ is the maximal level of $T$ containing some element of $X \cap Y$.
Then $t^+$ becomes the open ball centered at any fixed branch containing $t$ and with radius $1/n$, where $n$ is the height of $t$.
Note that with this metric, $\br T$ is always complete.
We will call $\br T$ the \emph{branch space} of the tree $T$.

\section{The Banach-Mazur game}

Historically, the game was invented by Mazur in 1935 (see Telg\'arsky~\cite{Telgarsky}) and can be described as follows.
Fix a set $X$ contained in the unit interval $\unii = [0,1]$.
Two players alternately choose non-degenerate intervals $J_0 \sups J_1 \sups \cdots$ contained in $\unii$.
The first player wins if $X \cap \bigcap_{\Ntr}J_n \nnempty$.
Otherwise, the second player wins.
Note that in Mazur's setting, the first player starts the game.
Thus, he has a winning strategy, for example, if the interior of $X$ is nonempty (actually, the minimal requirement is that $X$ is residual in some interval).
The situation changes drastically if the second player starts the game.
This was considered by Choquet in 1958, and the game was played in arbitrary topological spaces.
A theorem of Oxtoby (published in 1957) says that if the game is played in a metrizable space $X$ and the second player starts, then the first player has a winning strategy if and only if $X$ contains a dense completely metrizable subspace.

We propose a more general setting for the Banach-Mazur game.
First of all, note that the definition of the game requires only the structure $\pair{\Tau^+(X)}{\sups}$, that is, the result of a concrete play is an increasing sequence in this poset (or, a decreasing sequence when one considers $\subs$ instead of $\sups$).
Next, in order to say who wins, one needs to know which increasing sequences in $\pair{\Tau^+(X)}{\sups}$ are ``good'' for the first player.
This can be done by defining a bigger poset, where some of the sequences are bounded or have the least upper bound.

We now prepare the following framework for the Banach-Mazur game.
Namely, given a poset $\pair P\loe$, let $\ciagi P$ denote the poset of all countably generated ideals of $P$, ordered by inclusion.
Then $P$ can be naturally identified with a subset of $\ciagi P$, via the mapping
$$P \owns p \; \mapsto \; \icint p \in \ciagi P.$$
Before defining the Banach-Mazur game, in order to avoid confusion, we shall give names to the Players: \emph{Eve} and \emph{Odd}.
The result of a play will be a sequence $u_0 \loe u_1 \loe u_2 \loe u_3 \loe u_4 \loe \cdots$ in $P$, where the $u_n$ with $n$ even are chosen by Eve, while the $u_n$ with $n$ odd are chosen by Odd.
Below is the precise definition.

\begin{df}[The Banach-Mazur game in posets]
Fix a poset $\pair P \loe$ and fix $W \subs \ciagi P$.
The \emph{Banach-Mazur game} $\BMG P W$ is defined in the following way.
There are two players: \emph{Eve} and \emph{Odd}.
Eve starts the game by choosing $u_0 \in P$.
Odd responds by choosing $u_1 \in P$ with $u_0 \loe u_1$.
Then Eve responds by choosing $u_2 \in P$ with $u_1 \loe u_2$.
In general, if after one player's move we have a sequence $u_0 \loe \dots \loe u_n$ then the other player (no matter whether it is Eve or Odd) responds by choosing $u_{n+1} \in P$ with $u_n \loe u_{n+1}$.
We say that \emph{Odd wins} if the ideal generated by $\ciag u$ is an element of $W$; otherwise \emph{Eve wins}.
\end{df}

It is clear that if $\pair P \loe$ is of the form $\pair{\Tau^+(X)}{\sups}$, where $X$ is a topological space and $W$ consists of all ideals whose intersection is nonempty then the game defined above is just the classical Banach-Mazur game in the setting proposed by Choquet, where Eve is supposed to start the game.

Below we give an example from the theory of forcing.

\begin{ex}\label{EXfOrcing}
Let $\pair P \loe$ be a poset and let $\Dee$ be family of cofinal subsets of $P$.
An ideal $I$ of $P$ is \emph{$\Dee$-generic} if $I \cap D \nnempty$ for every $D \in \Dee$.
Let $W \subs \ciagi P$ be the family of all $\Dee$-generic ideals (a priori, we do not assume that $W \nnempty$).
We claim that if $\Dee$ is countable then Odd has a winning strategy in $\BMG P W$.

Indeed, let $\Dee = \ciag D$ and suppose $p_{2n}$ was the last Eve's choice in a fixed play.
Odd should choose $p_{2n+1} \in D_n$ so that $p_{2n+1} \goe p_{2n}$.
Applying this strategy, it is evident that Odd wins.
\end{ex}

The above example can be viewed as a strengthening of the well-known and simple Rasiowa-Sikorski lemma on the existence of generic ideals with respect to countably many cofinal sets.
Let us admit that in forcing theory typically the ordering is reversed and instead of a ``generic ideal" one uses the name ``generic filter".

\section{Rephrasing some classical results}\label{SecKlasik}

In this section we review some well-known results concerning the Banach-Mazur game in topological spaces, adapting them to our setting.

Given a poset $\pair P \loe$, we say that $A$ is an \emph{antichain} in $P$ if it consists of pairwise incompatible elements, where $x,y \in P$ are \emph{incompatible} if there is no $c \in P$ with $a \loe c$ and $b \loe c$; otherwise we say that $a$ and $b$ are \emph{compatible}.
An antichain $A$ is \emph{maximal} if it cannot be extended to a bigger antichain, that is, every element of $P \setminus A$ is compatible with some element of $A$.

\begin{tw}\label{ThmSDGsdbf}
Let $P$ be a poset, $W \subs \ciagi P$, and suppose that Odd has a winning strategy in $\BMG P W$.
Then there exists a metric tree $T \subs P$ with the following properties:
\begin{enumerate}
\item[(1)] Every level of $T$ is a maximal antichain in $P$.
\item[(2)] For every $I \in \br T$, the ideal generated by $I$ in $P$ is an element of $W$.
\end{enumerate}
\end{tw}

The proof of the theorem above is actually a direct translation of Oxtoby's arguments.
We shall see that under some circumstances the converse is also true.

\begin{pf}
Let $A_0$ be a maximal antichain in $P$ consisting of Odd's responses to the first Eve's move.
For each $a \in A_0$, choose a maximal antichain $A_1(a)$ in $\ciint a$ consisting of Odd's responses to the second Eve's move after $a$ (more formally, these are responses to $3$-element sequences where $a$ was the second element chosen by Eve).
We set $A_1 = \bigcup_{a \in A_0}A_1(a)$ and we note that $A_1$ is a maximal antichain in $P$.
Continuing this way, we obtain maximal antichains $\ciag A$, where $A_{n+1} = \bigcup_{a \in A_n}A_{n+1}(a)$ and $A_{n+1}(a)$ is a maximal antichain above $a$ consisting of Odd's responses to a suitable partial play.
Finally, $T = \bigcup_{\ntr}A_n \subs P$ is a metric tree satisfying (1).
Every branch $I$ of $T$ encodes a play of $\BMG P W$ where Odd was using his winning strategy, thus the ideal of $P$ generated by $I$ must be an element of $W$.
This completes the proof.
\end{pf}

\begin{tw}
Assume $P$, $W \subs \ciagi P$ and $T \subs P$ are as in Theorem~\ref{ThmSDGsdbf} above (in particular, $T$ is a metric tree satisfying (1), (2)).
If $W$ is a final segment in $\ciagi P$ then Odd has a winning strategy in $\BMG P W$.
\end{tw}

\begin{pf}
Let us describe the following strategy for Odd.
Assuming $u_0 \loe \dots u_{n}$ is a partial play with $n$ even, Odd chooses some element $a_{n/2}$ from the $(n/2)$th level of $T$ such that $a_{n/2} > a_{(n-1)/2}$ (in case $n>0$) and there is $v \in P$ satisfying $u_n \loe v$ and $a_{n/2} \loe v$.
He puts $u_{n+1} := v$.

After playing infinitely many steps of the game, we see that the ideal generated by the sequence $\ciag u$ contains the ideal generated by $\ciag a$ (which in turn is a cofinal subset of a branch of $T$), therefore it is in $W$, because $W$ is a final segment.
It follows that the strategy described above is winning for Odd.
\end{pf}

We say that a mapping of posets $\map \phi Q P$ is \emph{dominating} if
\begin{enumerate}
\item[(D1)] $\phi$ is order preserving, $\img \phi Q$ is cofinal in $P$, and
\item[(D2)] for every $q \in Q$, for every $p \in P$ with $\phi q \loe p$, there exists $q' \goe q$ in $Q$ such that $p \loe \phi q'$.
\end{enumerate} 

The following result allows us to ``move" the Banach-Mazur game from one poset to another, without changing its status.

\begin{tw}\label{Thmawfkbasfkj}
Let $\map \phi Q P$ be a dominating mapping of posets, let $W \subs \ciagi P$ and let $W^\phi$ consist of all ideals $I$ of $Q$ such that the ideal generated by $\img \phi I$ is in $W$.
The following conditions are equivalent:
\begin{enumerate}
\item[(a)] Odd has a winning strategy in $\BMG P W$.
\item[(b)] Odd has a winning strategy in $\BMG Q {W^\phi}$.
\end{enumerate}
The same applies to Eve.
\end{tw}

\begin{pf}
Suppose Odd has a winning strategy $\Sigma$ in $\BMG P W$.
We describe his winning strategy in $\BMG Q {W^\phi}$.
Namely, suppose Eve has chosen $v_0 \in Q$.
Odd first finds $u_1 \goe \phi v_0$ according to $\Sigma$ and then, using (D2), finds $v_1 \in Q$ such that $v_0 \loe v_1$ and $u_1 \loe \phi v_1$.
Finally, $v_1$ is Odd's response to the one-element sequence $v_0$.

In general, given a sequence $v_0 \loe \dots \loe v_{n-1}$, where $n>0$ is odd, we assume that we have the following sequence in $P$:
$$\phi v_0 \loe u_1 \loe \phi v_1 \loe \phi v_2 \loe u_3 \loe \phi v_3 \loe \phi v_4 \loe \dots \loe u_{n-2} \loe \phi v_{n-2} \loe \phi v_{n-1},$$
where $u_k = \Sigma(\phi v_0, u_1, \phi v_2, \dots, u_{k-2}, \phi v_{k-1})$ for every odd $k < n$.
Let
$$u_n = \Sigma(\phi v_0, u_1, \phi v_2, u_3, \phi v_4, \dots, u_{n-2}, \phi v_{n-1}).$$
Using (D2), Odd finds $v_n \goe v_{n-1}$ in $Q$ such that $u_n \loe \phi v_n$.
Finally, $v_n$ is Odd's response to $v_0 \loe \dots \loe v_{n-1}$.

Note that this strategy is winning because the $\phi$-image of the chain $v_0 \loe v_1 \loe \cdots$ is contained in a chain resulting from $\BMG P W$, where Odd was applying his winning strategy $\Sigma$.
In other words, the ideal generated by $\ciag {\phi v}$ is in $W$.
This shows the implication (a)$\implies$(b).

Suppose now that Odd has a winning strategy $\Pi$ in $\BMG Q {W^\phi}$.
We describe his winning strategy in $\BMG P W$.

Assume Eve has started with $u_0 \in P$.
Using (D1) Odd chooses $v_0 \in Q$ with $u_0 \loe \phi v_0$.
Next, he replies to $v_0$ according to $\Pi$, obtaining $v_1 \goe v_0$.
Then $u_1 := \phi v_1$ is his response to the one-element sequence $u_0$.

In general, given a sequence $u_0 \loe \dots \loe u_{n-1}$ with $n$ odd, we assume that there is a sequence $v_0 \loe \dots \loe v_{n-2}$ in $Q$ such that
$v_k = \Pi(v_0, \dots, v_{k-1})$ and $u_k = \phi(v_k) \loe u_{k+1}$ for every odd number $k < n$.
Odd's response to $u_0 \loe \dots \loe u_{n-1}$ is as follows.
Using (D2), he finds $v_{n-1} \goe v_{n-2}$ such that $\phi v_{n-1} \goe u_{n-1}$.
He takes $v_n = \Pi(v_0, \dots, v_{n-1})$ and responds with $u_n := \phi v_n$.

Note that the sequence $\ciag u$ resulting from this strategy contains a cofinal subsequence which is the $\phi$-image of $\ciag v$ which was winning in $\BMG Q {W^\phi}$), therefore the ideal generated by $\ciag u$ is in $W$.
We have shown the implication (b)$\implies$(a).

The second part (when Eve has a winning strategy) is almost the same, as the rules for both players are identical.
\end{pf}

\begin{wn}
Let $Q$ be a cofinal subset of a poset $P$ and let $W \subs \ciagi P$.
Then Odd / Eve has a winning strategy in $\BMG P W$ if and only if Odd / Eve has a winning strategy in $\BMG Q {W'}$, where $W' = \setof{I \cap Q}{I \in W}$.
\end{wn}

\begin{pf}
It suffices to notice that the identity mapping $\map \phi Q P$ is dominating.
\end{pf}

As a more concrete corollary, we see that Mazur was right by playing with nonempty open intervals instead of arbitrary open subsets of the real line.
Let us now recall Oxtoby's theorem \cite{Oxtoby}:

\begin{tw}\label{ThmOxtobybsdksd}
Let $X$ be a metrizable space. Then Odd has a winning strategy in $\BM X$ if and only if $X$ contains a dense completely metrizable subspace.
\end{tw}

\begin{pf}
Suppose first that $G \subs X$ is dense and completely metrizable and let $\rho$ be a complete metric on $G$.
We claim that Odd has a stationary winning strategy.
Namely, assuming $U$ was the last Eve's move, Odd responds with a nonempty open set $V$ satisfying the following two conditions: the closure of $V$ is contained in $U$, and the $\rho$-diameter of $V \cap G$ is finite, smaller than half of the $\rho$-diameter of $U \cap G$.
By Cantor's theorem, the intersection of any sequence resulting from a play with this strategy is a singleton of $G$.

Now suppose Odd has a winning strategy in $\BM X$ and let $\ovr X$ be the completion of $X$.
Define the following ordering on open sets: $U \rloe V$ iff either $U = V$ or else $\cl V \subs U$.
Note that the inclusion of $\pair {\Tau^(X)} \rloe$ in $\pair {\Tau^+(X)}{\sups}$ is dominating whenever $X$ is a regular space.
Let $Q = \pair {\Tau^+(\ovr X)} \rloe$, $P = \pair {\Tau^+(X)} \rloe$.
Let $\phi V = V \cap X$.
Then $\map \phi Q P$ is dominating. It is clear how to define $W \subs \ciagi P$ so that $\BM X$ becomes $\BMG P W$.
By Theorem~\ref{Thmawfkbasfkj}, we may consider $\BMG Q{W^\phi}$ instead.
Now let us look at Theorem~\ref{ThmSDGsdbf}.
Namely, we obtain a metric tree $T$ in $Q$, which translates to a tree of open sets in $\ovr X$ such that the intersection of each branch of $T$ is a single element of $X$ (because of Theorem~\ref{ThmSDGsdbf}(2)).
In other words, $T$ induces a dense completely metrizable subspace of $X$.
\end{pf}

We also have another variant of the Banach-Mazur game, for compact Hausdorff spaces.
Namely, if the classical Banach-Mazur game is played in a compact Hausdorff space, then Odd has an obvious stationary winning strategy: he always chooses an open set whose closure is contained in the last set chosen by Eve.
Now consider the Banach-Mazur game where the objective is to get a \emph{single} point in the intersection of the chain of open sets.
Let us call this game $\BMG X \star$, where $X$ is the topological space in question.
It turns out that there are non-metrizable compact Hausdorff spaces where Eve has a winning strategy in this game.

\begin{tw}
Let $K$ be a compact Hausdorff space.
The following properties are equivalent:
\begin{enumerate}
\item[(a)] Odd has a winning strategy in $\BMG K \star$.
\item[(b)] $K$ contains a dense $G_\delta$ metrizable subspace.
\end{enumerate}
\end{tw}

\begin{pf}
The proof of (b)$\implies$(a) is like the one in Theorem~\ref{ThmOxtobybsdksd}.
For the converse, we use the same tree $T$ as above, noting that it induces a metrizable subspace without assuming that the entire space $K$ is metrizable.
The fact that a compact Hausdorff space is regular is needed to conclude that the completely metrizable space of branches of $T$ is indeed dense in $K$.
\end{pf}

Let us recall the \emph{double arrow} space $K = D(\unii)$. This is a compact Hausdorff space whose universe is $((0,1] \times \sn 0) \cup ([0,1) \times \sn 1)$ endowed with the interval topology induced from the lexicographic ordering.
Let $\map p K \unii$ be the canonical projection.
Eve's winning strategy in $\BMG K \star$ is as follows: She always chooses an interval $U$ in $K$ such that $\img p U$ is in the interior of $\img p V$, where $V$ was the last choice of Odd.
Supposing that Odd wins while Eve plays this strategy, there would be a single point $x \in K$ in the intersection of $U_0 \sups U_1 \sups \cdots$ resulting from a play.
Now observe that $x$ is isolated from one side. For example, assume that $x = \pair y 0$, where $y<1$.
Then $\max \img p {U_n} = y$ from some point on (otherwise $\pair y1$ would be in the intersection), but this contradicts Eve's strategy saying that the closure of $\img p {U_{n+1}}$ is contained in $\img p {U_n}$.

At this point it is worth recalling that there are separable metric spaces in which the Banach-Mazur game is not determined.
Namely, recall that a \emph{Bernstein set} in a metrizable space is a set $S$ satisfying $S \cap P \nnempty \ne S \setminus P$ for every perfect set $P$ (a set is \emph{perfect} if it is nonempty, completely metrizable, and has no isolated points).
A Bernstein set in $2^\omega$ can be easily constructed by a transfinite induction, enumerating all perfect sets in $2^\omega$ and knowing that each perfect set in $2^\omega$ has cardinality continuum.

The following fact is well-known.

\begin{prop}
Let $X \subs 2^\omega$ be a Bernstein set.
Then the Banach-Mazur game $\BM X$ is not determined.
Namely, neither Eve nor Odd has a winning strategy in $\BM X$.
\end{prop}

\begin{pf}
By Theorem~\ref{Thmawfkbasfkj}, the game can be played in $\pair {2^\omega}\loe$, where $s \loe t$ means that $t$ extends $s$.
Furthermore, Odd wins if and only if the branch of $2^\omega$ resulting from a play corresponds to an element of $X$.
Suppose Odd has a winning strategy in $\BM X$.
Then there is a tree $T \subs 2^\omega$ satisfying conditions (1) and (2) of 
Theorem~\ref{ThmSDGsdbf}.
Note that the set of branches of $T$ is perfect, as all level of $T$ are maximal antichains in $2^\omega$.
This is a contradiction to the fact that $X$ is a Bernstein set.

Now suppose that Eve has a winning strategy in $\BM X$.
After her first move, we are in the same situation as before, interchanging the goals of the players.
Thus, a similar argument leads to a contradiction.
\end{pf}

\section{Applications to model theory}\label{SecFLimswegb}

We are now going to show that the Banach-Mazur game is determined when one considers the poset of all finitely generated structures of a fixed first order language, as long as some natural conditions are satisfied.

We now recall the concept of a \emph{\fra\ class}.
Namely, this a class $\fK$ of finitely generated models of a fixed language satisfying the following conditions:
\begin{enumerate}
\item[(F1)] For each $X, Y \in \fK$ there is $Z \in \fK$ such that both $X$ and $Y$ embed into $Z$.
\item[(F2)] Given embeddings $\map f Z X$, $\map g Z Y$ with $Z, X, Y \in \fK$, there exist $V \in \fK$ and embeddings $\map {f'}X V$, $\map {g'}Y V$ such that $f' \cmp f = g' \cmp g$.
\item[(F3)] For every $X \in \fK$, every finitely generated substructure of $X$ is in $\fK$.
\item[(F4)] There are countably many isomorphic types in $\fK$.
\end{enumerate}
Condition (F1) is called the \emph{joint embedding property}, (F2) is called the \emph{amalgamation property}.
Condition (F3) says that $\fK$ is \emph{hereditary} with respect to finitely generated substructures.

\fra\ theorem \cite{Fraisse} says that there exists a unique countably generated model $\U$ (called the \emph{\fra\ limit of $\fK$}) that can be presented as the union of a countable chain in $\fK$ and satisfies the following conditions:
\begin{enumerate}
\item[(U)] Every $X \in \fK$ embeds into $\U$.
\item[(E)] Given an isomorphism $\map h A B$ between finitely generated substructures $A, B \subs \U$, there exists an automorphism $\map H \U \U$ extending $h$.
\end{enumerate}

Let us denote by $\BMG \fK \U$ the Banach-Mazur game played in the poset\footnote{
We implicitly assume that all models ``live" in a certain fixed set, by this way we avoid dealing with a proper class instead of a set. For example, if finitely generated structures are finite then we may assume that the universe of each of them is a subset of $\Nat$.
}
consisting of all structures $X \in \fK$, where the ordering $\loe$ is inclusion (more precisely, extension of structures) and the winning ideals are precisely the structures isomorphic to the \fra\ limit of $\fK$.
In other words, Odd wins if and only if the resulting structure is isomorphic to $\U$.
In general, $\fK$ can be an arbitrary class and $\U$ can be an arbitrarily fixed model that is presentable as the union of a countable chain of models from the class $\fK$.
We shall later consider a more general version of this game, where a single model $\U$ is replaced by a family of models $\Yu$ and Odd wins if the union of the chain built by the two players is isomorphic to some $U \in \Yu$.
We denote this game by $\BMG \fK \Yu$.

Recall that a strategy of a fixed player is \emph{Markov} if his/her move depends only on the last move of the opponent and on the number of past moves.

\begin{tw}
If $\fK$ is a \fra\ class and $\U$ is the \fra\ limit of $\fK$, then
Odd has a Markov winning strategy in $\BMG \fK \U$.
\end{tw}

\begin{pf}
Let $\U$ be the \fra\ limit of $\fK$ and write $\U = \bigcup_{\ntr}U_n$, where $U_n \in \fK$ for each $\ntr$.
Odd's strategy is described as follows.

Supposing that last Eve's was $V_n$ (with $n$ even) and having recorded an embedding $\map{f_{n-1}}{V_{n-1}}{\U}$, Odd first chooses an embedding $\map g {V_n}{\U}$ extending $f_{n-1}$.
Next, he finds $V_{n+1} \in \fK$ with $V_n \loe V_{n+1}$ and an embedding $\map {f_{n+1}}{V_{n+1}}\U$ extending $g$ and such that $U_{n}$ is contained in the range of $f_{n+1}$.
In case $n=0$, we assume that $f_{-1}$ was the empty map.

It is clear that this strategy is winning for Odd, because after playing the game we obtain an isomorphism $f = \bigcup_{\ntr}f_{2n+1}$ of $\bigcup_{\ntr}V_n$ onto $\U$.
The strategy depends only on the result of last Eve's move and on the number of previous moves.
\end{pf}

One of the most important features of the \fra\ limit is that it is universal for the class of all countably generated structures obtained as unions of countable chains in $\fK$.
Using the Banach-Mazur game, we can give a simple direct argument in a more general setting.

\begin{tw}\label{Thmsdbgkjsege}
Let $\fK$ be a class of finitely generated models with the amalgamation property.
Let $\Yu$ be a class of countably generated models of the same language, such that Odd has a winning strategy in $\BMG \fK \Yu$.

Then every countably generated model representable as the union of a countable chain in $\fK$ is embeddable into some $U \in \Yu$.
\end{tw}

In case $\fK$ is a \fra\ class, we can set $\Yu = \sn \U$, where $\U$ is the \fra\ limit of $\fK$.

\begin{pf}
Assume $X = \bigcup_{\ntr}X_n$, where $X_n \loe X_{n+1}$ and $X_n \in \fK$ for each $\ntr$.
Let us play the game $\BMG \fK \Yu$, where Odd uses his winning strategy.
We shall denote by $U_0 \loe U_1 \loe U_2 \loe \cdots$ the concrete moves.
We describe a strategy of Eve leading to an embedding of $X$ into some $U \in \Yu$.

Namely, Eve starts with $U_0 = X_0$ and records $e_0 = \id{X_0}$.
Supposing that the last Odd's move was $U_{2n-1}$ and Eve has recorded an embedding $\map{e_{n-1}}{X_{n-1}}{U_{2n-2}}$, she uses the amalgamation property to find $U_{2n} \goe U_{2n-1}$ and an embedding $\map {e_n}{X_n}{U_{2n}}$ such that the diagram
$$\xymatrix{
U_{2n-2} \ar[r] & U_{2n-1} \ar[r] & U_{2n} \\
X_{n-1} \ar[u]^{e_{n-1}} \ar[rr] & & X_n \ar[u]_{e_n}
}$$
commutes, where the horizontal embeddings are inclusions.
By this way, $e_n$ extends $e_{n-1}$.

After playing the game, knowing that Odd wins, we conclude that $U = \bigcup_{\ntr}U_n$ is an element of $\Yu$ and $e = \bigcup_{\ntr}e_n$ is an embedding of $X$ into $U$.
\end{pf}

\separator

Let $\fK$ be a class of finitely generated models and let $U$ be a countably generated model of the same language such that Odd has a winning strategy in $\BMG \fK U$.
It is natural to ask what can be said about $\fK$ and $U$.

Clearly, $\fK$ has the joint-embedding property (F1), because Eve can play with any element of $\fK$, showing that $U$ is universal for $\fK$ (one can also use Theorem~\ref{Thmsdbgkjsege}, however this would be an overkill).
Assuming that $\fK$ consists of finite substructures, we conclude that (F4) must hold too, because $U$ is countable and therefore it has countably many isomorphic types of finite structures.
Obviously, (F3) may fail.
For example, let $\fG$ be a relational \fra\ class (say, the class of all finite graphs), and let $\fK$ be the subclass of $\fG$ consisting of all $G \in \fG$ whose cardinality is a prime number.
Let $\U$ be the \fra\ limit of $\fG$.
Then Odd has a winning strategy in the game $\BMG \fK \U$, as he can improve his Markov winning strategy for $\BMG \fG \U$ by enlarging his choices so that their cardinalities are always prime.

The following two examples from graph theory show that $\fK$ may fail the amalgamation property even when it satisfies (F1), (F3) and (F4).

\begin{ex}[Graphs with bounded degree]
Let $N>1$ be a fixed integer and let $\fK$ be the class of all finite graphs whose each vertex has degree $\loe N$.
It is well-known and easy to prove that each finite graph $G \in \fK$ embeds into a graph $H \in \fK$ such that the degree of \emph{every} vertex of $H$ is precisely $N$.
Let us call such a graph \emph{$N$-complete}.
For example, finite $2$-complete graphs are all cycles.

Clearly, $\fK$ is not a \fra\ class, as it fails the amalgamation property.

Let us enumerate by $\sett{H_n}{\ntr}$ all finite $N$-complete graphs.
Let $U = \bigoplus_{\ntr}U_n$, where ``$\bigoplus$" means the disjoint sum (and no extra edges between the summands), $\setof{U_n}{\ntr} = \setof{H_n}{\ntr}$, and for each $k \in \nat$ the set $\setof{\ntr}{U_n = H_k}$ is infinite.
In other words, $U$ is the direct sum of an indexed family consisting of (countably) infinitely many copies of each $H_n$.

We claim that Odd has a winning strategy in $\BMG \fK U$.
Indeed, after $n$th move of Eve resulting in a graph $G_{2n} \in \fK$, Odd chooses $G_{2n+1}$ of the form $\bigoplus_{i<k(n)}U_i$, knowing that each component of $G_{2n}$ embeds into some $U_i$.
He only to take care that $k(n) \to \infty$ while $n \to \infty$.
By this way, the graph resulting from a single play is obviously isomorphic to $U$.
\end{ex}

\begin{ex}[Cycle-free graphs]
Let $\fK$ denote the class of all finite cycle-free graphs.
Again, this is not a \fra\ class, as it fails the amalgamation property.
On the other hand, we claim that there is a countable cycle-free graph $B$ such that Odd has a winning strategy in $\BMG \fK B$.

Namely, let $T$ be the (uniquely determined) countable connected cycle-free graph whose each vertex has infinite degree.
The graph $T$ is well-known as the complete infinitely-branching tree with a single root.

Let $B$ be the direct sum of $\omega$ copies of $T$.
The winning strategy of Odd is as follows.
At stage $n$, after Eve's move $G_{2n}$, Odd chooses $G_{2n+1} \sups G_{2n}$ so that each component of $G_{2n+1}$ is a large enough part of $T$ (e.g. contains at least $n$ levels of $T$, when fixing the root).
Odd also takes care that at stage $n$ his graph $G_{2n+1}$ has at least $n$ components.
By this way Odd wins the play.
\end{ex}

\paragraph{Conclusion.}
Let, as above, $\fK$ be a countable class of finitely generated models of a fixed first order language and assume $U$ is such that Odd has a winning strategy in $\BMG \fK U$ (obviously, $U$ must be presentable as the union of a countable chain of models from $\fK$).
In that case we say that $U$ is \emph{generic} over $\fK$ (see Example~\ref{EXfOrcing} for an inspiration).
We have seen that \fra\ limits are generic over their \fra\ classes, however, there exist generic models that are not \fra\ limits in the usual sense.
There exists a category-theoretic generalization of \fra\ limits~\cite{K_Fraisse}, which in the case of models discards condition (F3) of being hereditary and possibly makes restrictions on embeddings.
By this way, we can talk about \emph{\fra\ categories} instead of \fra\ classes.
It can be proved that if $\fK$ contains a \fra\ subcategory $\fL$ that is \emph{dominating} in the sense of~\cite{K_Fraisse} then Odd has a winning strategy in $\BMG \fK U$ if and only if $U$ is the \fra\ limit of $\fL$ (in the setting of~\cite{K_Fraisse}).
We do not know whether the converse holds true.
In any case, generic objects seem to be a natural and applicable (see Theorem~\ref{Thmsdbgkjsege}) generalization of \fra\ limits.

Finally, let us note that in the Banach-Mazur game $\BMG \fK U$, the class $\fK$ can be just an abstract class of objects as long as the notion of an ``embedding" is defined.
It seems that the language of category theory is most suitable here.
Namely, $\fK$ could be a fixed category and $U$ could be a fixed object (typically in a bigger category containing $\fK$) that is isomorphic to the colimit of some sequence in $\fK$.
This approach will be explored elsewhere.

\end{document}